\documentclass[12pt,english]{smfart}

\DeclareMathAlphabet{\eusm}{U}{}{}{}  
\SetMathAlphabet\eusm{normal}{U}{eus}{m}{n}
\SetMathAlphabet\eusm{bold}{U}{eus}{b}{n}

\DeclareMathAlphabet{\eufrak}{U}{}{}{}  
\SetMathAlphabet\eufrak{normal}{U}{euf}{m}{n}
\SetMathAlphabet\eufrak{bold}{U}{euf}{b}{n}

\newtheorem{theorem}{Theorem}[section]
\newtheorem{proposition}[theorem]{Proposition}

\theoremstyle{definition}
\newtheorem{definition}[theorem]{Definition}
\newtheorem{example}[theorem]{Example}

\theoremstyle{remark}
\newtheorem{remark}[theorem]{Remark}
\numberwithin{equation}{section}

\newcommand{\bmultg}{\!\begin{array}{c} {\scriptstyle\times} \\[-12pt]\cup\end{array}\!}
\newcommand{\bmultk}{\!\!{\scriptstyle\begin{array}{c} {\scriptscriptstyle\times} \\[-12pt]\cup\end{array}}\!\!}

\author{Uwe Franz}
\address{D\'epartement de math\'ematiques de Besan\c{c}on,
Universit\'e de Franche-Comt\'e 16, route de Gray, 25030
Besan\c{c}on cedex, France}
\curraddr{Graduate School of Information
Sciences, Tohoku University, Sendai 980-8579, Japan}
\urladdr{http://www-math.univ-fcomte.fr/pp\underline{ }Annu/UFRANZ/}
\title[Multiplicative Boolean convolution]{Boolean convolution of probability measures on the unit circle}
\alttitle{Convolution bool\'eenne de probabilit\'es sur le cercle}

\begin{document}
\frontmatter

\begin{abstract}
We introduce the boolean convolution for probability measures on the unit circle. Roughly speaking, it describes the distribution of the product of two boolean independent unitary random variables. We find an analogue of the characteristic function and determine all infinitely divisible probability measures on the unit circle for the boolean convolution.
\end{abstract}

\begin{altabstract}
La convolution bool\'eenne de deux probabilit\'es sur le cercle est
d\'efinie comme la distribution du produit de deux op\'erateurs
unitaires $U$ et $V$ tels que $U-1$ et $V-1$. Un analogue de la
fonction caract\'eristique est donn\'ee et les lois infiniment
divisibles pour cette convolution sont caract\'eris\'ees.
\end{altabstract}

\subjclass{46L50; 60E05}

\keywords{Boolean independence, Boolean convolution, infinite
divisibility, L\'evy-Khintchine formula}

\altkeywords{Ind\'ependence bool\'eenne, convolution bool\'eenne,
lois infiniment divisibles, formule de L\'evy-Khintchine }

\thanks{Work supported in part by the European Community's Human Potential Programme under contract HPRN-CT-2002-00279 QP-Applications and a DAAD-KBN cooperation}

\maketitle

\mainmatter

\section{Introduction}

The boolean convolution of probability measures on the real line was introduced and studied in \cite{speicher+woroudi93}. In particular, it turned out that all probability measures on the real line are infinitely divisible for the boolean convolution. In this note we give a negative answer to Sch\"urmann's question, if the same is also true for the multiplicative boolean convolution of probability measures on the unit circle.

In Section \ref{sec-indep} we recall the definition of boolean independence and prove a crucial theorem on the distribution of two unitary operators $U$ and $V$, if $U-1$ and $V-1$ are boolean independent, see Theorem \ref{theo-F}.

The boolean convolution for probability measures on the unit circle is then defined in Section \ref{sec-conv}. We motivate the definition by the relation between the boolean product and the conditionally free product \cite{bozejko+leinert+speicher96}. But Definition \ref{def-bool-conv} also agrees with the theory of L\'evy processes on dual groups. I.e., if we take the algebra of polynomials on the unit circle as a dual group in the sense of Voiculescu \cite{voiculescu87} and define L\'evy processes on this dual group as in \cite{schuermann95b,franz03b}, then the distributions of the increments of these processes form convolution semigroups in the sense of Definition \ref{def-bool-conv}.

In Propositions \ref{prop-F} and \ref{prop-mult} we show that the transform $\mu\mapsto F_\mu$ with
\[
F_\mu(z)=\frac{1}{z}\,\frac{\psi_\mu(z)}{1+\psi_\mu(z)},
\]
where
\[
\psi_\mu(z)=\int_{S^1}\frac{xz}{1-xz}{\rm d}\mu(x),
\]
for $|z|<1$, reduces the multiplicative boolean convolution to the multiplication of holomorphic functions in the unit disk.

Finally, in Theorem \ref{theo-inf-div} we give a L\'evy-Khintchin formula that describes all probability measures on the unit circle that are infinitely divisible for the boolean convolution.

\section{Boolean independence}\label{sec-indep}

\begin{definition}
Let $H$ be a Hilbert space and $\varphi:\mathcal{B}(H)\to \mathbb{C}$ a normal state on $\mathcal{B}$. Two operators $X,Y\in\mathcal{B}(H)$ are called {\em boolean independent} w.r.t.\ $\varphi$, if
\begin{equation}\label{eq-def-boolean}
\varphi(X^{n_1} Y^{m_1} X^{n_2} \cdots X^{n_k} Y^{m_k})=\prod_{\ell=1}^k\varphi(X^{n_\ell})\varphi(Y^{m_\ell})
\end{equation}
holds for all $k\in \mathbb{N}$, $n_1,m_k\ge 0$, $n_2,\ldots,n_k,m_1,\ldots,m_{k-1}\ge 1$.
\end{definition}

For a unitary operator $U$ on $H$ we define functions $\psi_U$ and $F_U$ on the unit disk by
\[
\psi_U(z)=\varphi\left(\frac{z U}{1-zU}\right), \qquad |z|<1,
\]
and
\[
F_U(z)=\frac{1}{z}\,\frac{\psi_U(z)}{1+\psi_U(z)}, \qquad 0 < |z| < 1.
\]

The following theorem will be crucial for our study of the boolean convolution of probability measures in Section \ref{sec-conv}.

\begin{theorem}\label{theo-F}
Let $U,V$ be two unitary operators on $H$ such that $X=U-1$ and $Y=V-1$ are boolean independent w.r.t.\ $\varphi$. Then we have
\[
F_{UV}(z)=F_U(z)F_V(z)
\]
for all $z\in\mathbb{D}=\{z\in\mathbb{C}:|z|<1\}$.
\end{theorem}
\begin{proof}
Let $w\in\mathbb{C}$ with $|w|>4$.

Note that the resolvent $\eufrak{R}_{X+Y+XY}(w)=(w-X-Y-XY)^{-1}$ of $UV-1=X+Y+XY$ can be written as
\[
\eufrak{R}_{X+Y+XY}(w)=\eufrak{R}_X(w)\big(1-(X+1)Y\eufrak{R}_X(w)\big)^{-1}.
\]
We have chosen $|w|$ sufficiently large so that $\eufrak{R}_{X+Y+XY}(z)$ expands into the norm-convergent series
\[
\eufrak{R}_{X+Y+XY}(w) = \eufrak{R}_X(w)\sum_{n=0}^\infty \big((X+1)Y\eufrak{R}_X(w)\big)^n.
\]
Introducing the operators $\eufrak{K}_X(z)=\eufrak{R}_X(w)(X+1)-\frac{1}{w}$, we get
\begin{eqnarray*}
&& \eufrak{R}_{X+Y+XY}(w) \\
&=& \eufrak{R}_X(w)\left(1+(X+1)\sum_{n=1}^\infty \left(Y\left(\eufrak{K}_X(w)+\frac{1}{w}\right)\right)^{n-1}Y\eufrak{R}_X(w)\right) \\
&=& \eufrak{R}_X(w)\left(1+(X+1)\sum_{n=1}^\infty \sum_{{n_1+\cdots + n_k=n}\atop{n_1,\ldots,n_k\ge 1}} w^{-n+k}Y^{n_1}\eufrak{K}_X(w)Y^{n_2}\cdots Y^{n_k}\eufrak{R}_X(w)\right)
\end{eqnarray*}
The operator $\eufrak{K}_X(w)=\eufrak{R}_X(w)(X+1)-\frac{1}{w}$ can be written as a norm-convergent series
\[
\eufrak{K}_X(w)=\frac{X+1}{w-X}-\frac{1}{w}=\sum_{n=1}^\infty \left(\frac{1}{w^n}+\frac{1}{w^{n+1}}\right) X^n
\]
with vanishing constant term. Using Equation \eqref{eq-def-boolean}, we therefore have
\begin{gather*}
\varphi\big(\eufrak{R}_X(w)(X+1)Y^{n_1}\eufrak{K}_X(w)Y^{n_2}\cdots Y^{n_k}\eufrak{R}_X(w)\big) \\
= \left(K_X(w)+\frac{1}{w}\right) \varphi(Y^{n_1})\cdots \varphi(Y^{n_k}) \big(K_X(w)\big)^{k-1} G_X(w)
\end{gather*}
where $K_X(w)=\varphi\big(\eufrak{K}_X(w)\big)$ and $G_X(w)=\varphi\big(\eufrak{R}_X(w)\big)$.

We obtain the following formula for the Cauchy transform $G_{X+Y+XY}(w)=\varphi\big(\eufrak{R}_{X+Y+XY}(w)\big)$ of the distribution of $X+Y+XY$ in the state $\varphi$,
\begin{gather*}
G_{X+Y+XY}(w) \\
= G_X(w)\left(1+\left(K_X(w)+\frac{1}{w}\right)\sum_{n=1}^\infty \sum_{{n_1+\cdots + n_k=n}\atop{n_1,\ldots,n_k\ge 1}}\frac{\varphi(Y^{n_1})}{w^{n_1+1}}\cdots \frac{\varphi(Y^{n_k})}{w^{n_k+1}}w^{2k}\big(K_X(w)\big)^{k-1}\right)
\end{gather*}
Since $G_Y(w)-\frac{1}{w}=\sum_{n=1}^\infty \frac{\varphi(Y^n)}{w^{n+1}}$, we can rewrite this as
\begin{eqnarray*}
&& G_{X+Y+XY}(w) \\
&=& G_X(w)\left(1+\left(K_X(w)+\frac{1}{w}\right)\sum_{k=1}^\infty w^{2k}\left(G_Y(w)-\frac{1}{w}\right)^k\big(K_X(w)\big)^{k-1}\right) \\
 &=& G_X(w)\left(1+\frac{w^2\left(G_Y(w)-\frac{1}{w}\right)\left(K_X(w)+\frac{1}{w}\right)}{1-w^2\left(G_Y(w)-\frac{1}{w}\right)K_X(w)}\right)
\end{eqnarray*}
Substituting
\[
K_X(w)=(w+1)\left(G_X(w)-\frac{1}{w}\right)
\]
into this, we get after some simplification and a shift of the argument
\begin{eqnarray*}
G_{UV}(w) &=& G_{X+Y+XY}(w-1) \\
&=& \frac{G_U(w)G_V(w)}{w\big(G_U(w)+G_V(w)\big)-w(w-1)G_U(w)G_V(w)-1}
\end{eqnarray*}
for $|w|>4$. Substituting now
\[
G_Z(w)=\varphi\left(\frac{1}{w-Z}\right) = \frac{1}{w}\left(\psi_Z\left(\frac{1}{w}\right)+1\right)=\frac{1}{w-F_Z\left(\frac{1}{w}\right)}
\]
for $Z\in\{U,V,UV\}$ and replacing $w$ by $\frac{1}{z}$, this becomes finally
\[
F_{UV}(z) =F_U(z)F_V(z)
\]
for $0<|z|<\frac{1}{4}$. But by uniqueness of analytic continuation the identity holds for all $z\in\mathbb{D}$.
\end{proof}

\section{The Boolean convolution of measures on the unit circle}\label{sec-conv}

If $X\in\mathcal{B}(H)$ is self-adjoint, then there exists a uniquely determined compactly supported probability measure $\mu$ on $\mathbb{R}$ such that
\[
\varphi(X^k)=\int_\mathbb{R} x^k{\rm d}\mu
\]
holds for all $k\in\mathbb{N}$. We call this measure the distribution of $X$ in the state $\varphi$ and denote it by $\mathcal{L}_\varphi(X)=\mu$.

Let $X$ and $Y$ be two self-adjoint operators on $H$ and assume that $X$ and $Y$ are boolean independent w.r.t.\ $\varphi$. Then the joint moments $\varphi(X^{n_1} Y^{m_1} X^{n_2} \cdots X^{n_k} Y^{m_k})$ of the pair $X$ and $Y$ are uniquely determined by the distributions of $X$ and $Y$, since Equation \eqref{eq-def-boolean} can be used to compute them from the moments of $X$ and $Y$. In particular, the distribution of the sum $X+Y$ is uniquely determined by the marginal distributions $\mathcal{L}_\varphi(X)$ and $\mathcal{L}_\varphi(Y)$ and given by the boolean convolution
\[
\mathcal{L}_\varphi(X+Y)=\mathcal{L}_\varphi(X)\uplus\mathcal{L}_\varphi(Y).
\]
It was actually this property that motivated the definition of the boolean convolution in \cite{speicher+woroudi93}.

Let now $U$ be a unitary operator on $H$. Then there exists a unique probability measure $\mu$ on the unit circle such that
\[
\varphi(U^k) =\int_{S^1} x^k{\rm d}\mu(x)
\]
holds for all $k\in\mathbb{N}$. We will denote this measure by $\mathcal{L}_\varphi(U)=\mu$ and call it the distribution of $U$ in the state $\varphi$.

If $U,V$ are two boolean independent unitary operators, then the moments of their product are given by
\[
\varphi\big((UV)^n\big)=\varphi(UVUV\cdots UV)=\varphi(U)^n\varphi(V)^n, \qquad n\in\mathbb{N}.
\]
Since only the first moments of $U$ and $V$ enter, this does not lead to an interesting convolution of probability measures on the unit circle. The right definition can be deduced by recalling that the boolean product is a special case of the conditionally free product \cite{bozejko+leinert+speicher96}, which is defined for pairs of states. Looking at the algebra of continuous functions on the real line and taking evaluation at the origin for the second state, on recovers the additive boolean convolution. The important property of the choice of the second state is that the Dirac mass at the origin is an idempotent for the free convolution. If instead one takes the algebra of continuous functions on the unit circle and integration against the Dirac measure at $1$ for the second state, one arrives at the following definition of the multiplicative boolean convolution.

\begin{definition}\label{def-bool-conv}
Let $U,V\in \mathcal{B}(H)$ be two unitary operators with distributions $\mu=\mathcal{L}_\varphi(U)$, $\nu=\mathcal{L}_\varphi(V)$. Assume furthermore that $U-1$ and $V-1$ are boolean independent. Then we call the probability measure
\[
\mu\bmultg\nu=\mathcal{L}_\varphi(UV)
\]
the {\em boolean convolution} of $\mu$ and $\nu$. It is uniquely determined by $\mu$ and $\nu$.

We call a probability measure $\mu$ on the unit circle {\em infinitely divisible} w.r.t.\ the boolean convolution, if there exists a probability measure $\mu_n$ with
\begin{eqnarray*}
\mu &=& \underbrace{\mu_n\bmultg\mu_n\bmultg\cdots\bmultg\mu_n} \\
 & & \hspace*{10mm}n \mbox{ times}
\end{eqnarray*}
for all $n\in\mathbb{N}$.
\end{definition}

Let $\mu$ be a probability measure on the unit circle. Then it is uniquely characterized by the function
\[
\psi_\mu(z)=\int_{S^1} \frac{xz}{1-xz}{\rm d}\mu(x), \qquad |z|<1.
\]
We shall also need the function
\[
F_\mu(z)=\frac{1}{z}\,\frac{\psi_\mu(z)}{1+\psi_\mu(z)}, \qquad |z|<1.
\]
\begin{example}\label{exa-F}
\item[(a)]
Let $\mu=\delta_{e^{ib}}$ be a Dirac measure in $e^{ib}$, $0\le b<2\pi$. Then we get $\psi_{\delta_{e^{ib}}}(z)=\frac{ze^{ib}}{1-ze^{ib}}$ and $F_{\delta_{e^{ib}}}(z)=e^{ib}$ for $z\in\mathbb{D}$.
\item[(b)]
Let now $\mu=p\delta_{e^{ib_1}}+(1-p)\delta_{e^{ib_2}}$, with $0< p<1$ and $0\le b_1,b_2<2\pi$. Then we get
\[
\psi_\mu(z)=\frac{\left(pe^{ib_1}+(1-p)e^{ib_2}\right)z-e^{i(b_1+b_2)}z^2}{1-(e^{ib_1}+e^{ib_2})z+e^{i(b_1+b_2)}z^2},
\]
and
\[
F_\mu(z)=e^{i(b_1+b_2)}\frac{\alpha-z}{1-\overline{\alpha}z},
\]
for $z\in\mathbb{D}$, where $\alpha=pe^{-ib_2}+(1-p)e^{-ib_1}$.
\item[(c)]
Denote by $\lambda$ the uniform distribution on $S^1$, i.e.\ the Haar measure. Then we get $\psi_\lambda(z)=0$ and $F_\lambda(z)=0$ for $z\in\mathbb{D}$.
\item[(d)]
Consider now $\lambda^{(n)}=\frac{1}{n}\sum_{k=0}^{n-1}\delta_{z_k}$, where $z_k=e^{2\pi i k/n}$, for $n\in\mathbb{N}$. This is the Haar measure on the cyclic subgroup $\{z_k:k=0,\ldots,n-1\}$ of $S^1$. Here we get $\psi_{\lambda^{(n)}}=\frac{z^n}{1-z^n}$ and $F_{\lambda^{(n)}}(z)=z^{n-1}$ for $z\in\mathbb{D}$.
\end{example}

The following proposition shows that the map $\mu\mapsto F_\mu$ establishes an isomorphism between probability measures on the unit circle and holomorphic functions in the unit disk.

\begin{proposition}\label{prop-F}
The following are equivalent for a function $F$ on the unit disk $\mathbb{D}=\{z\in\mathbb{C}:|z|<1\}$.
\item[(i)]
We have $F:\mathbb{D}\to\overline{\mathbb{D}}$ and $F$ is holomorphic.
\item[(ii)]
There exists a probability measure $\mu$ on the unit circle $S^1=\{z\in\mathbb{C}:|z|=1\}$ such that
\[
F(z)=\frac{1}{z}\,\frac{\psi_\mu(z)}{1+\psi_\mu(z)}, \qquad z\in \mathbb{D},
\]
where
\[
\psi_\mu(z)=\int_{S^1}\frac{xz}{1-xz}{\rm d}\mu(x), \qquad z\in \mathbb{D}.
\]
The probability measure $\mu$ is uniquely determined by $F$.
\end{proposition}
\begin{remark}
\item[(a)]
If $F$ is not constant, then it is open and we can replace $\overline{\mathbb{D}}$ in part (i) by $\mathbb{D}$. If $F$ is constant and takes a value in $\partial \mathbb{D}=S^1$, then the associated measure $\mu$ is the Dirac measure in that point.
\item[(b)]
The function $\frac{\psi_\mu}{1+\psi_\mu}$ in (ii) has a zero at the origin and therefore $F=\frac{1}{z}\,\frac{\psi_\mu}{1+\psi_\mu}$ on $\mathbb{D}\backslash\{0\}$ has a unique analytic extension to $\mathbb{D}$
\end{remark}
\begin{proof}
(i)$\Rightarrow$(ii): Let $F:\mathbb{D}\to\overline{\mathbb{D}}$ be holomorphic. then
\[
g(z)=2 \frac{zF(z)}{1-zF(z)} +1, \qquad z\in \mathbb{D},
\]
is also holomorphic and has non-negative real part. Therefore $g$ has a Herglotz representation
\begin{equation}\label{eq-nu}
g(z)=ib+\int_{S^1}\frac{w+z}{w-z}{\rm d}\nu(w), \qquad z\in \mathbb{D},
\end{equation}
see \cite[Theorem V.4.10]{bhatia97}. Here $b={\rm Im}\, g(0)=0$ and $\nu(S^1)={\rm Re}\, g(0)=1$, i.e.\ $\nu$ is probability measure on the unit circle. This implies
\[
\psi(z) = \frac{1}{2}\big(g(z)-1\big) =
\int_{S^1}\frac{z/w}{1-z/w}{\rm d}\nu(w)=
\int_{S^1}\frac{xz}{1-xz}{\rm d}\mu(x)
\]
where $\mu=\tilde{\nu}$ is the image of $\nu$ under the map $x\mapsto x^{-1}$.

The uniqueness of $\mu$ follows from the uniqueness of the Herglotz representation.

(ii)$\Rightarrow$(i): The M\"obius transformations $\sigma_x(z)=\frac{xz}{1-xz}$ map the unit disk to the half plane with real part bigger $-\frac{1}{2}$ for all $x\in S^1$. Therefore $\psi(z)=\int_{S^1}\frac{xz}{1-xz}{\rm d}\mu(x)$ is holomorphic and satisfies ${\rm Re}\,\psi_\mu(z)>-\frac{1}{2}$. This implies that $h=\sigma_1^{-1}\circ\psi$ is holomorphic and maps the unit disk to itself. Since $\psi(0)=0$, we also have $h(0)=0$ and $F(z)=\frac{1}{z}h(z)$, $0<|z|<1$, can be extended analytically to the origin. The function $F$ satisfies
\[
|F(z)|=\left|\frac{h(z)}{z}\right| \le (1+\varepsilon) \quad\mbox{ for }\quad |z|=\frac{1}{1+\varepsilon},\quad \varepsilon>0.
\]
Letting $\varepsilon$ go to zero, one gets $|F(z)|\le 1$ for all $|z|<1$ by the maximum principle.
\end{proof}

We can now restate Theorem \ref{theo-F} for the boolean convolution.

\begin{proposition}\label{prop-mult}
Let $\mu,\nu$ be probability measures on the unit circle. Then we have
\[
F_{\mu\bmultk\nu}(z)=F_\mu(z)F_\nu(z)
\]
for $|z|<1$.
\end{proposition}
\begin{proof}
Follows immediately from Theorem \ref{theo-F} and Definition \ref{def-bool-conv}.
\end{proof}

The Haar measure $\lambda$ has $F_\lambda\equiv 0$. It is infinitely divisible for the boolean convolution, since $\lambda^{\bmultk n}=\lambda$ for all $n\in\mathbb{N}$, and satisfies $\lambda\bmultg\mu=\lambda$ for all probability measures $\mu$ on the unit circle.

The following theorem describes all other infinitely divisible probability measures by a L\'evy-Khintchin-type formula.

\begin{theorem}\label{theo-inf-div}
A probability measure $\mu\not=\lambda$ on the unit circle is infinitely divisible for the boolean convolution if and only if there exists a real number $b\in[0,2\pi)$ and a finite measure $\rho$ on $S^1$ such that
\[
F_\mu(z)=\exp\big(u(z)\big)
\]
with
\[
u(z)=ib-\int_{S^1}\frac{x+z}{x-z}{\rm d}\rho(x)
\]
for $|z|<1$.

The pair $(b,\rho)$ is called the {\em characteristic pair} of $\mu$, it is uniquely determined by $\mu$.
\end{theorem}
\begin{proof}
A probability measure $\mu$ on the unit circle is infinitely divisible for the boolean convolution if and only if there exist holomorphic functions $F_n$ with
\[
\big(F_n(z)\big)^n =F_\mu(z), \qquad z\in\mathbb{D}
\]
for all $n\in\mathbb{N}$. If $F$ is not identically zero, then this happens if and only if $0\not\in F(\mathbb{D})$. In this case we can take the logarithm of $F$, i.e.\ there exists a holomorphic function $u$ with ${\rm Re}\,u(z)\le 0$ for $z\in\mathbb{D}$ such that
\[
F(z)=\exp\big(u(z)\big),\qquad z\in\mathbb{D}.
\]
This function $u$ is uniquely determined up to an integer multiple of $2\pi i$. We choose $u$ such that $0\le{\rm Im}\,u(0)<2\pi$. Since the real part of $-u$ is non-negative, it has a Herglotz representation
\[
u(z)=ib-\int_{S^1}\frac{x+z}{x-z}{\rm d}\rho(x), \qquad z\in\mathbb{D}.
\]
with $b={\rm Im}\,u(0)\in[0,2\pi)$ and $\rho$ a finite measure on the unit circle.

Conversely, given a pair $(b,\rho)$ of a real number $b\in[0,2\pi)$ and $\rho$ a finite measure on $S^1$, the functions
\[
F(z)=\exp\big(u(z)\big) \quad \mbox{ and }\quad F_n(z)=\exp\big(u(z)/n\big), \qquad z\in\mathbb{D},
\]
with
\[
u(z)=ib-\int_{S^1}\frac{x+z}{x-z}{\rm d}\rho(x), \qquad z\in\mathbb{D}.
\]
give by Proposition \ref{prop-F} probability measures $\mu$ and $\mu_n$ which satisfy
\[
\mu=\mu_n^{\bmultk n}
\]
by Proposition \ref{prop-mult}. Therefore $\mu$ is infinitely divisible for the boolean convolution.

The uniqueness follows from the uniqueness of the Herglotz representation.
\end{proof}
\begin{remark}
\item[(a)]
It follows that not all probability measures on the unit circle are infinitely divisible for the boolean convolution. Let $\mu\not=\lambda$. If the first moment of $\mu$ vanishes, then $\psi_\mu$ and $\frac{\psi_\mu}{1+\psi_\mu}$ have a zero of order at least two at the origin. Therefore $F_\mu(0)=0$, but $F_\mu\not\equiv0$, and we can not take arbitrary holomorphic roots of $F_\mu$. This is equivalent to $\mu$ not being infinitely divisible.
\item[(b)]
It follows that every infinitely divisible probability measure $\mu$ can be embedded into a continuous convolution semigroup, just take for $(\mu_t)_{t\ge0}$ the convolution semigroup with the measures $\mu_t$ given by
\[
F_{\mu_t}(z)=\exp\left( itb - t\int_{S^1}\frac{x+z}{x-z}{\rm d}\rho(x)\right),\qquad |z|<1,
\]
for $t\ge0$, if $(b,\rho)$ is the characteristic pair of $\mu$. The embedding is not unique, because replacing $b$ by $b'=b+2\pi k$, $k\in\mathbb{Z}\backslash\{0\}$, leads to another convolution semigroup $(\mu'_t)_{t\ge 0}$ which also satisfies $\mu'_1=\mu$.
\item[(c)]
Probability measures with $\mu\bmultg\mu=\mu$ are called idempotents. For the classical convolution on the unit circle these are exactly the Haar measures of compact subgroups of $S^1$, i.e.\ the probability measures given in Example \ref{exa-F} (d). Idempotents for the boolean convolution are characterized by the condition $\big(F_\mu(z))^2=F_\mu(z)$ for $|z|<1$. Since $F_\mu$ has to be holomorphic, this leaves only the constant functions $F\equiv 1$ and $F\equiv 0$, which correspond to the Dirac measure at $1$ and the Haar measure on $S^1$.
\item[(d)]
Analytic functions on the unit disk have a unique factorization into a product $F=BSO$, where $B$ is a Blaschke product, $S$ is a singular function, i.e.\ of the form
\[
S(z)=\exp\left(-\int_{S^1}\frac{w+z}{w-z}{\rm d}\tau(w)\right),\qquad z\in\mathbb{D},
\]
with a measure $\tau$ that is singular w.r.t.\ to the Haar measure, and $O$ is an outer function, i.e.\ of the form
\[
O(z)=c\exp\left(-\int_{S^1}\frac{w+z}{w-z}q(w){\rm d}\lambda(w)\right),\qquad z\in\mathbb{D},
\]
with $c\in\mathbb{C}$, $|c|=1$, and $q$ a real-valued integrable function, see, e.g., \cite{hoffman62}. If $|F(z)\le 1$ for $|z|<1$, then $q$ is non-negative. The function $B$ has the form
\[
B(z)=z^p\prod_{n=1}^\infty\left(\frac{\overline{\alpha}_n}{|\alpha_n|}\,\frac{\alpha_n-z}{\,1-\overline{\alpha}_nz\,}\right)^{p_n}, \qquad z\in\mathbb{D},
\]
where $p,p_1,p_2,\ldots$ are non-negative integers and $(\alpha_n)_{n\in\mathbb{N}}$ is a sequence such that $\alpha_n\in\mathbb{D}\backslash\{0\}$, $\alpha_n\not=\alpha_m$ for $n\not=m$, and $\prod_{n=1}^\infty |\alpha_n|^{p_n}$ converges. By Proposition \ref{prop-mult} and Theorem \ref{theo-inf-div}, we get a corresponding factorization
\[
\mu=\mu_B\bmultg\mu_S\bmultg\mu_O
\]
of any probability measure $\mu$ on the unit circle into the boolean convolution of a not-infinitely divisible measure $\mu_B$, an infinitely divisible singular measure $\mu_S$ and an infinitely divisible measure $\mu_O$ that is absolutely continuous w.r.t.\ $\lambda$.
\end{remark}
\begin{example}
\item[(a)]
For $\rho=0$ we get the Dirac measure $\delta_{e^{ib}}$.
\item[(b)]
Let $b=0$ and $\rho=\delta_x$, $|x|=1$. Then we get $F_\mu(z)=\exp\left(\frac{z+x}{z-x}\right)$ and $\mu$ is characterized by
\[
\int_{S^1} \frac{1+wz}{1-wz}{\rm d}\mu=1+2\psi_\mu(z)=\frac{1+zF_\mu(z)}{1-zF_\mu(z)}=\frac{\exp\left(\frac{x+z}{x-z}\right)+z}{\exp\left(\frac{x+z}{x-z}\right)-z}=:k(z),
\]
for $|z|<1$. Since the real part of $k(z)$ goes to zero almost everywhere when $z$ goes to the unit circle, the measure $\mu$ is singular with respect to the Haar measure $\lambda$. Denote by $\{z_n:n\in\mathbb{N}\}$ the zeros of $\exp\left(\frac{x+z}{x-z}\right)-z$ on the unit circle. One can show that they have an accumulation point at $x$. Since $g$ has single poles in each $z_n$, the measure $\mu$ has a Dirac mass $a_n\delta_{\overline{z}_n}$ at $\overline{z}_n$. Let $x=e^{i\beta}$, $z_n=e^{i\beta_n}$, then the coefficient $a_n$ is given by
\begin{eqnarray*}
a_n &=& \lim_{{z\to z_n}\atop{z\in\mathbb{D}}}\frac{z_n-z}{z_n+z}k(z)=\lim_{{z\to z_n}\atop{z\in\mathbb{D}}}\frac{z_n-z}{\exp\left(\frac{x+z}{x-z}\right)-z} \\
&=& \lim_{{z\to z_n}\atop{z\in\mathbb{D}}}\frac{-1}{\frac{2x}{(x-z)^2}\exp\left(\frac{x+z}{x-z}\right)-1} = \frac{(x-z_n)^2}{x^2+z_n^2-4xz_n} \\
&=& \frac{1-\cos(\beta-\beta_n)}{2-\cos(\beta-\beta_n)}.
\end{eqnarray*}
We do not know, if $\mu$ also charges the point $\overline{x}$.
\item[(c)]
Consider now the constant function $F\equiv re^{ib}$ with $0\le r<1$, $0\le b<2\pi$. This function also belongs to an infinitely divisible measure. It is absolutely continuous with respect to the Haar measure $\lambda$, with density given by the Poisson kernel
\[
P_r(\beta)=\frac{1}{2\pi}{\rm Re}\,\frac{1+\overline{z}re^{i b}}{1-\overline{z}re^{ib}} =
\frac{1}{2\pi}\frac{1-r^2}{1+r^2-2r\cos(\beta-b)},
\]
in $z=e^{i\beta}\in S^1$. The measure $\rho$ in the characteristic pair of $\mu$ is determined by the condition
\[
\int_{S^1}\frac{x+z}{x-z}{\rm d}\rho =-\ln r.
\]
It is equal to $\rho=-\ln r \lambda$, i.e.\ it is a scalar multiple of the Haar measure.
\end{example}

\section*{Acknowledgment}

I wish to thank Rolf Gohm for his remarks and advice.

\backmatter

\end{document}